\documentclass{amsart}
\usepackage{amsmath,amssymb,amsthm,amscd}
\newtheorem{formula}{}[section]
\newtheorem{proposition}[formula]{Proposition}
\newtheorem{corollary}[formula]{Corollary}
\newtheorem{lemma}[formula]{Lemma}
\newtheorem{theorem}[formula]{Theorem}
\theoremstyle{definition}
\newtheorem{definition}[formula]{Definition}
\newtheorem{example}[formula]{Example}
\theoremstyle{remark}
\newtheorem*{remark}{Remark}

\begin{document}

\title{Graded filiform Lie algebras and symplectic nilmanifolds}
\author{Dmitri V. Millionschikov}
\thanks{Partially supported by
the Russian Foundation for Fundamental Research, grant no. 99-01-00090 and PAI-RUSSIE, dossier no. 04495UL}
\subjclass{17B30, 17B56, 17B70, 53D}
\address{Department of Mathematics and Mechanics, Moscow
State University, 119899 Moscow, RUSSIA}
\curraddr{Universit\'e Louis Pasteur, UFR de Math\'ematique et d'Informatique, 7 rue Ren\'e Descartes - 67084 Strasbourg Cedex (France)}
\email{million@mech.math.msu.su}

\begin{abstract}
We study symplectic (contact) structures on nilmanifolds that correspond to the 
filiform Lie algebras - nilpotent Lie algebras
of the maximal length of the descending central sequence. We give a complete
classification of filiform Lie algebras that possess a basis $e_1, \dots, e_n$,
$[e_i,e_j]=c_{ij}e_{i{+}j}$ ($\mathbb N$-graded Lie algebras).
In particular we describe the spaces of symplectic cohomology classes for all even-dimensional algebras of the list.
It is proved that a symplectic filiform Lie algebra $\mathfrak g$ is a filtered deformation of
some $\mathbb N$-graded symplectic filiform Lie algebra ${\mathfrak g}_0$. 
But this condition is not sufficient. 
A spectral sequence
is constructed in order to answer the question whether a given deformation 
of a $\mathbb N$-graded symplectic filiform Lie algebra ${\mathfrak g}_0$ admits a symplectic structure or not.  
Other applications and examples are discussed.
\end{abstract}
\date{}

\maketitle

\section*{Introduction}
Nilmanifolds $M=G/\Gamma$ (compact homogeneous spaces of nilpotent Lie groups $G$
over lattices $\Gamma$) are widely used in
topology and geometry. We can mention for instance Gromov's theory of Carnot-Carath\'eodory spaces, 
the examples of Anosov flows, as well as different
counter-examples in symplectic and complex topology.
The study of invariant geometric structures (Riemannian metrics, symplectic (contact) structures,
complex structures) on nilmanifolds reduces to the problems in terms of the tangent nilpotent  Lie algebra $\mathfrak{g}$ of $G$.
The first examples of compact  
symplectic manifolds with no K\"ahler structure were nilmanifolds
(see \cite{TO}, \cite{CFG} for references).
I. Babenko and I. Taimanov considered  in \cite{BT1}, \cite{BT2} an interesting family of nilmanifolds $M_n$
that they used for construction of 
examples of non formal simply connected symplectic manifolds.
In their proof I. Babenko and I. Taimanov used the existence of integer symplectic structure on
$M_n$ and a non-trivial triple Massey product in $H^1(M_n, \mathbb R)$. 
The problem of computation of corresponding $H^*(M_n, \mathbb R)$
was partially solved in \cite{Mill1}, \cite{Mill2}. Using  Nomizu's theorem ~\cite{Nz} the bigraded
structure of $H^*(M_n)=H^*(\mathcal{V}_n)$ was studied.  
$\mathcal{V}_n$ denotes the corresponding nilpotent Lie algebra.

The natural idea was to find other similar examples of nilpotent symplectic Lie algebras. 
But the theory of nilponet Lie algebras has no effective classification 
tools like in the theory of semi-simple Lie algebras. 
For instance we have Morosov's classification \cite{Mor}
of nilpotent Lie algebras with dimensions $\le 6$. Morosov's table of Lie algebras was used 
in \cite{GozeB}, \cite{GozeKh} 
for studing symplectic structures and complex structures in \cite{Sal}.
Articles on classification
of $7$-dimensional complex nilpotent Lie algebras regularily appear and however we have no certitude 
that the problem is finally solved. We restrict ourselves to the class of so-called filiform Lie algebras --
nilpotent Lie algebras $\mathfrak{g}$ with a maximal length $s{=}\dim \mathfrak{g}{-}1$ of the descending central sequence of
$\mathfrak{g}$. The Lie algebras
$\mathcal{V}_n$ are filiform ones.
The study of this class of nilpotent Lie algebras was started by M.~Vergne in 
\cite{V1}, \cite{V2} and later was continued by M.~Goze and his coauthors 
(see \cite{GozeKh}, \cite{AG}, \cite{GozeB}, \cite{GJKh1}, \cite{GJKh2}, \cite{Kh1}).
But in these articles  filiform Lie algebras were studied over $\mathbb C$.
In particular in \cite{GJKh2} symplectic filiform Lie algebras were classified in dimensions $\le 10$ 
but again over complex field $\mathbb C$. 

In the present article we mainly work with $\mathbb K = \mathbb R$ but it is possible to consider
an arbitrary field $\mathbb K$ of zero characteristic as it was done in \cite{Br}.

This paper is organized as follows. In the sections 1 -- 4 we review all necessary facts of nilpotent Lie
algebras theory as well as the theory of nilmanifolds. Namely we discuss
two filtrations of a filiform Lie algebra $\mathfrak{g}$: the canonical one $C$,
defined by the ideals of the descending central sequence and the second one $L$, defined 
by means of so-called adapted basis $e_1, \dots, e_n$ of $\mathfrak{g}$. In the sequel we consider
two corresponding associated graded Lie algebras ${\rm gr}_C\mathfrak{g}$ and ${\rm gr}_L\mathfrak{g}$ that
play a very important role.

We start our study of filiform Lie algebras in the Section 5 with the special case of $\mathbb N$-graded filiform Lie algebras,
i.e. $\mathbb N$-graded Lie algebras $\mathfrak{g}=\oplus_{\alpha}\mathfrak{g}_{\alpha}$ such that 
$\dim \mathfrak{g}_i=1,  \; 1 \le i \le n$,
$\mathfrak{g}_i =0, \; i > n$,
$[\mathfrak{g}_1,\mathfrak{g}_i]=\mathfrak{g}_{i{+}1}, \; i \ge 2$. We give a complete
classification (Theorem \ref{osnovn}) of these algebras.
 We recall that it was Y.~Khakimdjanov who discovered in 
\cite{Kh1} that there exists only a finite number of non-isomorphic $\mathbb N$-graded filiform Lie
algebras over $\mathbb C$ in dimensions $\ge 12$. We add in dimensions
$\ge 12$ three new classes $\mathfrak{m}_{0,1}(2k{+}1)$, $\mathfrak{m}_{0,2}(2k{+}2)$, $\mathfrak{m}_{0,3}(2k{+}3)$
that were missed in \cite{Kh1}, \cite{{GozeKh}} as well as explicit formulae for one-parameter families in
dimensions $\le 11$ (for an arbitrary field $\mathbb K$ of zero characteristic). In the proof we 
use an inductive procedure of construction of filiform Lie algebras by means of one-dimensional central extensions.    

In the Section 6, 7 we study symplectic structures on filiform Lie algebras. We describe all symplectic
$\mathbb N$-graded filiform Lie algebras (Theorem \ref{list_sympl}). Namely for a $\mathbb N$-graded filiform Lie algebra
$\mathfrak{g}$ we essentially use the bigraded structure of $H^*(\mathfrak{g})$.
We prove that if $\mathfrak{g}$ is symplectic than ${\rm gr}_C\mathfrak{g} \cong \mathfrak{m}_0(2l)$
and hence the filtration $L$ is defined. Let us consider $\mathbb N$-graded filiform Lie algebra
${\rm gr}_L\mathfrak{g}$. Theorem \ref{sympl_def} states that if  $\mathfrak{g}$ is symplectic than
the $\mathbb N$-graded algebra ${\rm gr}_L\mathfrak{g}$ is symplectic also and hence  
${\rm gr}_L\mathfrak{g}$ is isomorphic to one of the following algebras 
$$\mathfrak{m}_0(2l), \;{\mathcal V}_{2k}, \; {\mathfrak g}_{8, \alpha}, \; {\mathfrak g}_{10, \alpha},$$
where by ${\mathfrak g}_{8, \alpha}, \; {\mathfrak g}_{10, \alpha}$ we denote the one-parameter families of filiform Lie algebras 
in dimensions $8$ and $10$ respectively. But the last condition is not sufficient, an example
of a deformation of the symplectic $\mathbb N$-graded algebra $\mathfrak{m}_0(2l)$ that does not admit any symplectic structure
is discussed. A criterion of existence of symplectic structure in terms of the spectral sequence $E$ that corresponds to the
filtration $L$ is formulated (Theorem \ref{criterion}).

In the last Section we discuss other possible applications: contact structures on odd-dimensional filiform Lie algebras,
the classification problem of filiform Lie algebras and G.~Lauret's \cite{Lr} construction of nilmanifolds 
admitting Anosov diffeomorphism.  
   
\section{Symplectic and contact nilmanifolds}

\begin{definition}  A nilmanifold $M$ is a compact homogeneous space
of the form
$G/\Gamma,$ where $G$ is a simply connected nilpotent Lie group
and $\Gamma$ is a lattice in $G$.
\end{definition}

\begin{example}
A $n-$dimensional torus $T^n={\mathbb R}^n/{\mathbb Z}^n$.
\end{example}
\begin{example}
The Heisenberg manifold $M_3={\mathcal H}_3/\Gamma_3,$ where
${\mathcal H}_3$ is the group of all matrices of the form
$$
   \left( \begin{array}{lcr}
   1 & x & z\\
   0 & 1 & y\\
   0 & 0 & 1\\
   \end{array} \right) , ~~~ x,y,z \in {\mathbb R},
$$
and a lattice $\Gamma_3$ is a subgroup of matrices with $x,y,z
\in {\mathbb Z}$.
\end{example}

\begin{example}
We define the generalized Heisenberg group $H(1,p)$ as the group of all matrices of the form
$$
   \left( \begin{array}{lcr}
   I_p & X & Z\\
   0 & 1 & y\\
   0 & 0 & 1\\
   \end{array} \right) , 
$$
where $I_p$ denotes the identity $p \times p$ matrix, $X,Z$ are $p \times 1$ matrices of real numbers,
and $y$ is a real number. 
Let us consider a product $G(p,q)= H(1,p) \times H(1,q)$. By $\Gamma(p,q)$ we denote a lattice of all
matrices with integer entries. We can define a family nilmanifolds $M(p,q)=G(p,q)/\Gamma(p,q)$,
where $\dim M(p,q)=2p+2q+2$.
\end{example}

\begin{theorem}[A.I. Malcev \cite{Mal}]
Let $G$ be a simply connected nilpotent Lie group with a tangent Lie
algebra $\mathfrak{g}$.
Then $G$ has a co-compact lattice $\Gamma$ (i.e. $G/\Gamma$ is a compact
space) if and only if there exists a 
basis $e_1, e_2, \dots, e_n$ in $\mathfrak{g}$
such that the constants $\{c_{ij}^{k} \}$ of Lie structure
$[e_i,e_j]= c_{ij}^{k}e_k$
are all rational numbers.
\end{theorem}

\begin{example}
\label{malcev}
Let $\mathfrak{g}$ be a real nilpotent Lie algebra
with a basis $e_1, e_2, \dots, e_n$ and rational constants of Lie structure:
$[e_i,e_j]=c_{ij}^k e_k, c_{ij}^k \in \mathbb Q$.
One can also define 
a group structure $*$ in the vector space $\mathfrak{g}$
using the Campbell-Hausdorff formula:

$$x*y=x+y+\frac{1}{2}[x,y]+\frac{1}{12}([x,[x,y]]-[y,[y,x]])+\dots$$

The nilpotent group $G=(\mathfrak{g},*)$
has a lattice $\Gamma$ such that 
the quotient space $M=G/\Gamma$ is
a compact manifold.
\end{example} 

\begin{definition}  An even-dimensional manifold $M^{2k}$ is called symplectic if it has
a closed $2$-form $\omega$ 
(symplectic structure on $M$) such that $\omega$ is non-degenerate $\forall P \in M^{2k}$.
\end{definition}
One of the first non-trivial examples of symplectic nilmanifolds was
Kodaira (or Kodaira-Thurston) surface $K=M_3 \times S^1$ with a symplectic form 
$\omega =\frac{1}{2 \pi}dy \wedge d \phi + dx \wedge dz$, where by $M_3$ we denote
the Heisenberg manifold and by $d\phi$ the standard volume form on $S^1$.
The symplectic manifold $K=M_3 \times S^1$ does not admit a K\"ahler structure (see \cite{TO} for references).

\begin{theorem}[C.~Benson, C.~Gordon \cite{BG}]
If a nilmanifold $G/\Gamma$ admit a K\"ahler structure, then $G$ is abelian and
$G/\Gamma$ is diffeomorphic to a torus. 
\end{theorem}

\begin{lemma}[\cite{CFG}]
The $2$-form $F$ on $G(p,q)$ defined by
$$F=dX^t_1 \wedge (dZ_1-X_1dy_1)+dX^t_2 \wedge (dZ_2-X_2dy_2)+dy_1 \wedge dy_2$$
defines a symplectic structure on $G(p,q)$ and $M(p,q)$, where
by $X_i, Z_i, y_i, i=1,2$ we denote global coordinates on the first $H(1,p)$ and the second
factor $H(1,q)$ of $G(p,q)$ respectively.
\end{lemma}
The family $M(p,q)$ was probably the first example of a infinite family of (non-formal) symplectic nilmanifolds (\cite{CFG}).

\begin{example}
Let ${\mathcal V}_n$ be a nilpotent Lie algebra
with a basis $e_1, e_2, \dots, e_n$ and a Lie bracket:
$$[e_i,e_j]= \left\{\begin{array}{r}
   (j-i)e_{i{+}j},  i+j \le n;\\
   0, i+j > n.\\
   \end{array} \right. $$
Let us denote by $M_n$ the corresponding nilmanifold that we can obtain
using Malcev's method \ref{malcev}.
\end{example} 

\begin{lemma}[\cite{BT1}]
$M_{2k}$ is a symplectic nilmanifold with a left-invariant symplectic
form $$\omega_{2k{+}1}=(2k{-}1)e^1 \land e^{2k}+
(2k{-}3)e^2 \land e^{2k{-}1}+
\dots{+}e^k \land e^{k{+}1},$$
where $e^1, \dots, e^{2k}$
is a dual basis of left-invariant 1-forms to the $e_1, \dots, e_{2k}.$
\end{lemma}

I. Babenko,
I. Taimanov in \cite{BT1}, \cite{BT2} applied the symplectic
blow-up procedure
to the symplectic embedding $M_{2k} \to {\mathbb C}P^{2k{+}1}$
in order to construct examples of simply connected non-formal symplectic manifolds.

We can identify the deRham complex $\Lambda^*(G/\Gamma)$
with subcomplex
$\Lambda^*_{\Gamma}(G) \subset \Lambda^*(G)$
of left-invariant forms on $G$ with respect to the action of the lattice
$\Gamma$. $\Lambda^*_{\Gamma}(G)$ has a
subcomplex $\Lambda^*_{G}(G)$ of left-invariant forms 
with respect to the whole action of $G$.
$\Lambda^*_{G}(G)$ is naturally isomorphic
to the Lie algebra cochain complex
$\Lambda^*(\mathfrak{g})$.
One can consider the corresponding inclusion
$$\psi: \Lambda^*(\mathfrak{g}) \to 
\Lambda^*(G/\Gamma).$$

\begin{theorem}[Nomizu ~\cite{Nz}]
Let $G/\Gamma$ be a compact nilmanifold then
the inclusion $\psi: \Lambda^*(\mathfrak{g}) \to 
\Lambda^*(G/\Gamma)$ induces the isomorphism 
$\psi^*: H^*(\mathfrak{g}) \to H^*(G/\Gamma, {\mathbb R})$
in cohomology.
\end{theorem}

\begin{definition}
An skew-symmetric non-degenerate bilinear form $\omega$ on the $2k$-dimensional Lie algebra $\mathfrak{g}$
is called symplectic if it closed: $d\omega=0$.
\end{definition}

It follows that $\omega^k =  \underbrace{\omega \wedge \omega \wedge \dots \wedge \omega}_{k} \ne 0$. 

\begin{corollary}
Let $G/\Gamma$ be a compact nilmanifold with a symplectic form $\omega$ then
$\omega = \psi(\tilde \omega)+ d \xi,$ where $\tilde \omega $ is a symplectic form on 
the tangent Lie algebra $\mathfrak{g}$ of $G$
and $\xi$ is a $1$-form on $G/\Gamma$.
\end{corollary}

Lie groups with left-invariant symplectic, K\"ahler and Poisson structures were discussed
in \cite{LM}. B.-Y. Chu proved in \cite{Chu} that a unimodular symplectic Lie algebra has to be solvable.  

\begin{definition}  An odd-dimensional manifold $N^{2k{+}1}$ is called contact if it has
a $1$-form $\alpha$ (the contact form) such that $\alpha \wedge d \alpha^k$ is a volume form on 
$N^{2k{+}1}$.
\end{definition}

\begin{theorem}[Theorem 3.2. in \cite{TO}]
\label{TO3.2}
Every symplectic nilmanifold ($M, \omega$) arises  as the orbit space of a free circle action on
a contact nilmanifold ($N, \alpha$) such that $\alpha$ is invariant under the action
and $\pi^*(\omega)=d \alpha$, where $\pi: N \to M$ is the orbit space projection.
\end{theorem}

\begin{theorem}[Theorem 3.3. in \cite{TO}]
\label{TO3.3}
Every contact nilmanifold with invariant contact form has a free circle
action whose orbit space is a symplectic manifold.
\end{theorem}

\section{$\mathbb N$-graded Lie algebras}
\begin{definition}
A Lie algebra $\mathfrak{g}$ is called $A$-graded,
if it is decomposed in a direct sum of subspaces such that
$$\mathfrak{g}=\bigoplus_{\alpha} \mathfrak{g}_{\alpha}, \; \alpha \in A,
\quad \quad [\mathfrak{g}_{\alpha}, \mathfrak{g}_{\beta}]
\subset \mathfrak{g}_{\alpha+\beta}, \;
\forall \: \alpha, \beta \in A.
$$
\end{definition}
Ž

Usually the following sets $A$ are considered: 
abelian groups ${\mathbb Z}$
and ${\mathbb Z}_m$.
 In this article we'll treat of $A={\mathbb N}$ -- abelian semigroup of natural numbers.

It is evident that $\mathbb N$-graded Lie algebras of finite dimension are nilpotent ones.

In the present article we will consider 
${\mathbb N}$-graded Lie algebras $\mathfrak{g}$ of the following types 

1) finite-dimensional Lie algebras $\mathfrak{g}=\oplus_{\alpha=1}^n {\mathfrak{g}}_{\alpha}$
$$\dim \mathfrak{g}_{\alpha}= \left\{\begin{array}{r}
   1, \; \alpha \le n;\\
   0, \; \alpha > n.\\
   \end{array} \right.; $$

2) infinite-dimensional Lie algebras $\mathfrak{g}=\oplus_{\alpha=1} {\mathfrak{g}}_{\alpha}$
$$\dim \mathfrak{g}_{\alpha}= 1, \; \alpha \ge 1. $$
 
\begin{example}
€The Lie algebra $\mathfrak{m}_0$ is defined 
by its infinite basis $e_1, e_2, \dots, e_n, \dots $
with commutating relations:
$$ [e_1,e_i]=e_{i+1}, \; \forall \; i \ge 2.$$
\end{example}

\begin{remark}
We will omit in the sequel trivial commutating relations $[e_i,e_j]=0$ in the definitions of Lie algebras.  
\end{remark}

\begin{example}
€The Lie algebra $\mathfrak{m}_2$ is defined by its infinite basis $e_1, e_2, \dots, e_n, \dots $
and commutating relations:
$$
[e_1, e_i ]=e_{i+1}, \quad \forall \; i \ge 2; \quad \quad
[e_2, e_j ]=e_{j+2}, \quad \forall \; j \ge 3. 
$$
\end{example}

\begin{example}
Let us define the algebra $L_k$ as the Lie algebra of polynomial vector fields
on the real line ${\mathbb R}^1$ with a zero in $x=0$ of order not less then $k+1$.
\end{example}
€The algebra   $L_k$ can be defined by its basis and  commutating relations
$$e_i=x^{i+1}\frac{d}{dx}, \; i \in {\mathbb N},\; i \ge k, \quad \quad
[e_i,e_j]= (j-i)e_{i{+}j}, \; \forall \;i,j \in {\mathbb N}.$$

€The algebras $\mathfrak{m}_0, \mathfrak{m}_2, L_1$ are
${\mathbb N}$-graded Lie algebras generated by two elements $e_1, e_2$.

€In her article A.Fialowski ~\cite{Fial} classified all  
${\mathbb N}$-graded Lie algebras of the type 2) with two generators.
In particular the following theorem holds on.
ˆ\begin{theorem}[€A. Fialowski, ~\cite{Fial}] 
\label{fialowski}
Let $\mathfrak{g}=\bigoplus_{\alpha=1}^{\infty}\mathfrak{g}_{\alpha}$
be a ${\mathbb N}$-graded Lie algebra such that:
\begin{equation}
\dim \mathfrak{g}_{\alpha}= 1, \; \alpha \ge 1; \quad
[\mathfrak{g}_{1}, \mathfrak{g}_{\alpha}]=\mathfrak{g}_{\alpha+1},
\; \forall \alpha \ge 2.
\end{equation}
'Then $\mathfrak{g}$ is isomorphic to one (and only one) Lie algebra from three given ones:
$$\mathfrak{m}_0, \; \mathfrak{m}_2, \; L_1.$$ 
\end{theorem}

\begin{remark}
Later some of Fialowski's results were rediscovered in \cite{Kh2}.
\end{remark}

Let  $\mathfrak{g}$ be a infinite-dimensional ${\mathbb N}$-graded Lie algebra,
then the quotient algebra $\mathfrak{g}(n) =
\mathfrak{g} / \oplus_{\alpha=n+1}^{\infty} \mathfrak{g}_{\alpha}$
is the $\mathbb N$-graded $n$-dimensional Lie algebra and  
$\mathfrak{g}= \lim_n \mathfrak{g}(n)$.

So we have three sequences of finite-dimensional $\mathbb N$-graded
Lie algebras: $\mathfrak{m}_0(n)$, $\mathfrak{m}_2(n)$, ${\mathcal V}_n$.
We took the notation 
${\mathcal V}_n=L_1/L_{n{+}1}$ as  
$L_1(n)$ has another meaning in the theory of infinite-dimensional Lie algebras.

\section{Nilpotent and filiform Lie algebras}

The sequence of ideals of a Lie algebra $\mathfrak{g}$
$$C^1\mathfrak{g}=\mathfrak{g} \; \supset \;
C^2\mathfrak{g}=[\mathfrak{g},\mathfrak{g}] \; \supset \; \dots
\; \supset \;
C^{k}\mathfrak{g}=[\mathfrak{g},C^{k-1}\mathfrak{g}] \; \supset
\; \dots$$
is called the descending central sequence of $\mathfrak{g}$.

A Lie algebra $\mathfrak{g}$ is called nilpotent if
there exists $s$ such that:
$$C^{s+1}\mathfrak{g}=[\mathfrak{g}, C^{s}\mathfrak{g}]=0,
\quad C^{s}\mathfrak{g} \: \ne 0.$$
'The natural number $s$ is called the nil-index of the nilpotent Lie algebra $\mathfrak{g}$, or 
$\mathfrak{g}$ is  called  $s$-step nilpotent Lie algebra.

Let $\mathfrak{g}$ be a Lie algebra. We call a set $F$ of subspaces
$$\mathfrak{g} \supset \dots \supset F^{i}\mathfrak{g} \supset
F^{i{+}1}\mathfrak{g} \supset \dots \qquad (i \in \mathbb Z)
$$
a decreasing filtration $F$ of $\mathfrak{g}$ if $F$ is compatible with the Lie structure
$$[F^{k}\mathfrak{g},F^{l}\mathfrak{g}] \subset F^{k+l}\mathfrak{g},\; \forall k,l \in \mathbb Z.$$

Let $\mathfrak{g}$ be a filtered Lie algebra.  
A graded Lie algebra 
$${\rm gr}_F\mathfrak{g}=\bigoplus_{k=1} ({\rm gr}_F\mathfrak{g})_k, \;\;
({\rm gr}_F\mathfrak{g})_k=F^{k}\mathfrak{g}/F^{k{+}1}\mathfrak{g}$$
is called the associated (to $\mathfrak{g}$) graded Lie algebra
${\rm gr}_F\mathfrak{g}$.

ˆ\begin{example}
The ideals $C^{k}\mathfrak{g}$ of the descending central sequence define
a decreasing filtration $C$ of the Lie algebra
$\mathfrak{g}$ 
$$C^{1}\mathfrak{g}=\mathfrak{g} \supset C^{2}\mathfrak{g} \supset
\dots \supset C^{k}\mathfrak{g} \supset \dots; \qquad
[C^{k}\mathfrak{g},C^{l}\mathfrak{g}] \subset C^{k+l}\mathfrak{g}.$$
"One can consider the associated graded Lie algebra
${\rm gr}_C\mathfrak{g}$.
\end{example}

The finite filtration $C$ of a nilpotent Lie algebra $\mathfrak g$ is called 
the canonical filtration of
a nilpotent Lie algebra $\mathfrak g$. 

\begin{proposition}
Let $\mathfrak{g}$ be a $n$-dimensional nilpotent Lie algebra.
Then for its nil-index we have the estimate $s \le n-1$.
\end{proposition}

"It follows from another estimate for nilpotent Lie algebras:
$$\dim ({\rm gr}_C\mathfrak{g})_1=
\dim (\mathfrak{g}/[\mathfrak{g},\mathfrak{g}]) \: \ge \: 2. $$

\begin{definition}
A nilpotent $n$-dimensional Lie algebra $\mathfrak{g}$ is  called
filiform Lie algebra if it has the nil-index $s=n-1$.  
\end{definition}

The Lie algebras $\mathfrak{m}_0(n)$, $\mathfrak{m}_2(n)$, ${\mathcal V}_n$
considered above are filiform Lie algebras.

' \begin{lemma}[ ŒM. Vergne \cite{V1}]
Let $\mathfrak{g}$ be a filiform Lie algebra, then there exists a basis
$e_1,e_2, \dots, e_n$ such that: 
\begin{equation}
\label{cycl_basis}
ad(e_1)(e_i)=[e_1,e_i]=e_{i+1}, \quad i=2,\dots,n-1.
\end{equation}
\end{lemma}

ŠCertainly a filiform Lie algebra $\mathfrak{g}$ besides (\ref{cycl_basis}) may have other non-trivial commutating relations. 
'In this sense the filiform Lie algebra  $\mathfrak{m}_0(n)$ is the simplest possible
example of a filiform Lie algebra.

\begin{proposition}
Let $\mathfrak{g}$ be a filiform Lie algebra and  ${\rm gr}_C\mathfrak{g}= \oplus_i ({\rm gr}_C\mathfrak{g})_i$ is
the corresponding associated (with respect to the  canonical filtration $C$)
graded Lie algebra. Then 
$$
\dim ({\rm gr}_C\mathfrak{g})_1=2, \quad
\dim ({\rm gr}_C\mathfrak{g})_2=\dots=
\dim ({\rm gr}_C\mathfrak{g})_{n{-}1}=1.
$$
\end{proposition}

ˆWe have the following isomorphisms of graded Lie algebras:
$$ {\rm gr}_C\mathfrak{m}_2(n) \cong
{\rm gr}_C{\mathcal V}_n \cong {\rm gr}_C\mathfrak{m}_0(n) \cong
\mathfrak{m}_0(n).$$

\begin{theorem}[ŒM. Vergne \cite{V2}]
\label{V_ne1}
Let $\mathfrak{g}=\oplus_{\alpha}\mathfrak{g}_{\alpha}$ be a
graded $n$-dimensional filiform Lie algebra and 
\begin{equation}
\label{can_grad}
\dim \mathfrak{g}_1=2, \quad
\dim \mathfrak{g}_2=\dots=
\dim \mathfrak{g}_{n{-}1}=1.
\end{equation}
'then  

1) if $n=2k+1$, then $\mathfrak{g}$ is isomorphic to 
$\mathfrak{m}_0(2k+1)$;

2) if $n=2k$, then $\mathfrak{g}$ is isomorphic either to
$\mathfrak{m}_0(2k)$ or to the Lie algebra $\mathfrak{m}_1(2k)$,
defined by its basis $e_1, \dots, e_{2k}$ 
and commutating relations:
$$
[e_1, e_i ]=e_{i+1}, \; i=2, \dots, 2k{-}1; \quad \quad
[e_j, e_{2k{+}1{-}j} ]=(-1)^{j{+}1}e_{2k}, \quad j=2, \dots, k.
$$
\end{theorem}

\begin{remark}
In the settings of the Theorem \ref{V_ne1}
the gradings of the algebras $\mathfrak{m}_0(n)$, $\mathfrak{m}_1(n)$
are defined as
$\mathfrak{g}_1=Span( e_1, e_2 ) $,
$\mathfrak{g}_i=Span( e_{i{+}1} ), \: i=2, \dots, n{-}1.$
\end{remark}

\begin{corollary}[ŒM. Vergne \cite{V2}]
Let $\mathfrak g$ be a filiform Lie algebra. Then one can choose a so-called
adapted basis $e_1, e_2, \dots, e_n$ in $\mathfrak g$:
\begin{equation}
[e_i, e_j]= \left\{\begin{array}{r}
   \sum \limits_{k{=}0}^{n{-}i{-}j}c_{ij}^k e_{i{+}j{+}k}, \; i+j \le n;\\
   (-1)^i \alpha e_n, \; i+j=n+1 ; \\
   0, \hspace{2.76em} \; i+j > n+1.
   \end{array} \right. 
\end{equation}
where $\alpha =0$ if $n$ is odd number.
\end{corollary}

\begin{example}
Let $\mathfrak g$ be a filiform Lie algebra     
with a fixed adapted basis $e_1, e_2, \dots, e_n$ such that $\alpha=0$. 
Then one can define a new (non-canonical) filtration $L$ of $\mathfrak g$:
'$$\mathfrak g=L^1 \mathfrak g \supset L^2 \mathfrak g \supset \dots \supset L^n \mathfrak g \supset \{0\},$$
where 
$L^k\mathfrak{g}{=}Span ( e_{k}, \dots, e_n ),
\, k=1,{\dots},n.$ 
\end{example}

\begin{remark}
The adapted basis is completely determined by $e_1, e_2$:
$e_3=[e_1, e_2], e_4=[e_1, e_3], \dots, e_n=[e_1, e_{n{-}1}]$. 
It is convenient do define filiform Lie algebras by means of their adapted basises.
\end{remark}

\section{Lie algebra cohomology}

Let us consider the cochain complex of a Lie algebra $\mathfrak{g}$ with 
$\dim \mathfrak{g}=n$:
$$
\begin{CD}
\mathbb K @>{d_0{=}0}>> \mathfrak{g}^* @>{d_1}>> \Lambda^2 (\mathfrak{g}^*) @>{d_2}>>
\dots @>{d_{n-1}}>>\Lambda^{n} (\mathfrak{g}^*) @>>> 0
\end{CD}
$$

where $d_1: \mathfrak{g}^* \rightarrow \Lambda^2 (\mathfrak{g}^*)$
is a dual mapping to the Lie bracket
$[ \, , ]: \Lambda^2 \mathfrak{g} \to \mathfrak{g}$,  
and the differential $d$ (that is a collection of $d_p$)
is a derivation of the exterior algebra $\Lambda^*(\mathfrak{g}^*)$,
that continues $d_1$:

$$
d(\rho \wedge \eta)=d\rho \wedge \eta+(-1)^{deg\rho} \rho \wedge d\eta,
\; \forall \rho, \eta \in \Lambda^{*} (\mathfrak{g}^*).
$$
"The condition $d^2=0$ is equivalent to the Jacobi identity in $\mathfrak{g}$.

ŠThe cohomology of $(\Lambda^{*} (\mathfrak{g}^*), d)$ is called the cohomology (with trivial coefficients) of the Lie algebra 
 $\mathfrak{g}$ and is denoted by $H^*(\mathfrak{g})$.

'The exterior algebra  $\Lambda^* \mathfrak{g}$ of a graded Lie algebra  
$\mathfrak{g}=\oplus_{\alpha}\mathfrak{g}_{\alpha}$
can be equipped with the second grading 
$\Lambda^* \mathfrak{g} =
\bigoplus_{\mu} \Lambda^*_{(\mu)} \mathfrak{g}$, where
$\Lambda^p_{(\mu)} \mathfrak{g}$ is spanned by 
$\xi_1 {\wedge} \xi_2 {\wedge} \dots {\wedge}
\xi_p \in  \Lambda^p \mathfrak{g}$ such that 
$$\xi_1 \in \mathfrak{g}_{\alpha_1}, \xi_2 \in \mathfrak{g}_{\alpha_2},
\dots,  \xi_p \in \mathfrak{g}_{\alpha_p}, \;
\alpha_1{+}\alpha_2{+}\dots{+}\alpha_p=\mu .$$

ŽLet us identify $\Lambda^* (\mathfrak{g}^*)$ with the space
$(\Lambda^* \mathfrak{g})^*$. Then  
$\Lambda^* (\mathfrak{g}^*)$ has also second grading 
$\Lambda^* (\mathfrak{g}^*) =
\bigoplus_{\lambda} \Lambda^*_{(\lambda)} (\mathfrak{g}^*)$,
where the subspace 
$\Lambda^p_{(\lambda)} (\mathfrak{g}^*)$ is defined by
$$ \Lambda^p_{(\lambda)} (\mathfrak{g}^*)= \left\{ \omega
\in \Lambda^p (\mathfrak{g}^*)
 \; | \; \omega(v)=0, \: \forall v \in
 \Lambda^p_{(\mu)} (\mathfrak{g}), \: \mu \ne \lambda \right\}.$$

The second grading is compatible with the differential  $d$ and with the exterior product
$$d \Lambda^p_{(\lambda)} (\mathfrak{g}^*)
\subset \Lambda^{p{+}1}_{(\lambda)} (\mathfrak{g}^*), \qquad
\Lambda^{p}_{(\lambda)} (\mathfrak{g}^*) \cdot
\Lambda^{q}_{(\mu)} (\mathfrak{g}^*) \subset
\Lambda^{p{+}q}_{(\lambda{+}\mu)} (\mathfrak{g}^*)
$$
'The exterior product in $\Lambda^*(\mathfrak{g}^*)$ 
induces the structure of bigraded algebra in cohomology $H^*(\mathfrak{g})$:
$$
H^{p}_{(\lambda)} (\mathfrak{g}) \otimes
H^{q}_{(\mu)} (\mathfrak{g}) \to
H^{p{+}q}_{(\lambda{+}\mu)} (\mathfrak{g}).
$$

\begin{example}
Let $\mathfrak{g}$ be a Lie algebra with the basis
 $e_1, e_2, \dots, e_n$ and commutating relations
$$[e_i,e_j]= \left\{\begin{array}{r}
   c_{ij}e_{i{+}j},  i+j \le n;\\
   0, i+j > n.\\[3pt]
   \end{array} \right. $$

Let us consider the dual basis $e^1, e^2, \dots, e^n$.
'We have the following formulae for the differential $d$:  
$$de^k=\sum_{i{+}j{=}k} c_{ij}e^i \wedge e^{j}, \; k=1,\dots,n.$$  
Now we can introduce the grading (that we will call the weight)
of $\Lambda^*(\mathfrak{g}^*)$:
$$\Lambda^* (\mathfrak{g}^*)=
\bigoplus_{\lambda{=}1}^{n(n{+}1)/2}
\Lambda^*_{(\lambda)} (\mathfrak{g}^*),$$
where a subspace 
$\Lambda^{p}_{(\lambda)} (\mathfrak{g}^*)$ is spanned by $p$-forms
$\{ e^{i_1} {\wedge} \dots {\wedge} e^{i_p}, \;
i_1{+}\dots{+}i_p {=} \lambda \}$. 
For instance a monomial 
$e^{i_1} \wedge \dots \wedge e^{i_p}$ has the degree $p$ and the weight 
$\lambda=i_1{+}\dots{+}i_p$.
\end{example}

Now we consider a filtered Lie algebra $\mathfrak{g}$ with 
a finite decreasing filtration $F$.
One can define a decreasing filtration $\tilde F$ of $\Lambda^* ({\mathfrak{g}}^*)$.  
$$ \tilde F^{\mu} \Lambda^p (\mathfrak{g}^*)= \left\{ \omega
\in \Lambda^p (\mathfrak{g}^*)
 \; | \; \omega(F^{\alpha_1} \mathfrak{g} \wedge \dots \wedge F^{\alpha_p}\mathfrak{g})=0, \:  
 \alpha_1{+} \dots {+} \alpha_p {+} \mu \ge 0 \right\}.$$

\begin{example}
Let $\mathfrak{g}$ be a Lie algebra with the basis
 $e_1, e_2, \dots, e_n$ and commutating relations
$$[e_i,e_j]= \left\{\begin{array}{r}
   \sum\limits_{k=0}^{n{-}i{-}j} c^k_{ij}e_{i{+}j{+}k},  i+j \le n;\\
   0, i+j > n.\\[3pt]
   \end{array} \right. $$

As it was remarked above the corresponding filtration $L$ can be defined.
The associated graded Lie algebra ${\rm gr}_L \mathfrak g$ can be defined by structure relations
$$[e_i,e_j]= \left\{\begin{array}{r}
   c_{ij}^0 e_{i{+}j},  i+j \le n;\\
   0, i+j > n.\\[3pt]
   \end{array} \right. $$
Let us consider the dual basis $e^1, e^2, \dots, e^n$. Then 
'$\tilde L^{\mu}\Lambda^{p} (\mathfrak{g}^*)$ is spanned by $p$-monomials of weights
less or equal to $-\mu$, i.e. by
$p$-forms $e^{i_1} {\wedge} \dots {\wedge} e^{i_p}$ such that $i_1{+}\dots{+}i_p {\le} -\mu$. For instance
'$$\tilde L^{{-}5}\Lambda^{2} (\mathfrak{g}^*)=Span(e^1 \wedge e^2, e^1 \wedge e^3, e^1 \wedge e^4, e^2 \wedge e^3).$$  
\end{example}

\begin{remark}
One can consider the spectral sequence $E_r$ that corresponds to the filtration $\tilde F$ of
the complex $\Lambda^* ({\mathfrak{g}}^*)$. We have an isomorphism (see \cite{V2} for example)
$$E^{p,q}_1 = H^{p{+}q}_{(-p)}({\rm gr}_F \mathfrak g).$$  
\end{remark}

\begin{definition}
€A Lie algebra $\mathfrak{g}$ is called unimodular
if the trace of a operator  $ad(x)$ is equal to zero for all $x \in \mathfrak{g}$.
\end{definition}
As an elemetary corollary of the Engel theorem we get that 
a nilpotent Lie algebra $\mathfrak{g}$ is unimodular.
\begin{theorem}[J.-L. Koszul \cite{K}]
Let $\mathfrak{g}$ be a $n$-dimensional unimodular  Lie algebra then
'\begin{enumerate}
\item $H^n(\mathfrak{g}) \cong \Lambda^n(\mathfrak{g}) \cong
Span( e^1 \wedge e^2 \wedge \dots \wedge e^n ) \cong {\mathbb K}$;
\item the multiplication $H^p(\mathfrak{g}) \otimes H^{n{-}p}(\mathfrak{g})
\to H^n(\mathfrak{g}) \cong {\mathbb K}$ defines a non-degenerate bilinear form
(Poincar\'e duality).
\end{enumerate}
\end{theorem}

\section{Filiform Lie algebras as one-dimensional central extensions}

Recall that a one-dimensional central extension of a Lie algebra $\mathfrak{g}$ is an exact sequence
\begin{equation}
\label{exactseq}
0 \to \mathbb K \to \tilde {\mathfrak g}  \to \mathfrak{g} \to 0 
\end{equation}
of Lie algebras and their homomorphisms, in which the image of the homomorphism
${\mathbb K} \to \tilde {\mathfrak{g}}$ 
is contained in the center of the Lie algebra $\mathfrak{g}$.
To the cocycle $c \in \Lambda^2(\mathfrak{g}^*)$ corresponds the extension
$$
0 \to {\mathbb K} \to {\mathbb K} \oplus \mathfrak{g} \to \mathfrak{g} \to 0
$$
where the Lie bracket in $\mathbb{K} \oplus \mathfrak{g}$ is defined by the formula
$$
[(\lambda,g), (\mu, h)]=(c(g,h), [g,h]).
$$
It can be checked directly that the Jacobi identity for this Lie bracket is equivalent to $c$
being cocycle and that to cohomologous cocycles correspond equivalent (in a obvious sense)
extensions.

If $\mathfrak{g}$ is a Lie algebra of the finite dimension $n$ defined by its basis $e_1, \dots, e_n$ and
structure relations
$$ [e_i,e_j]=\sum_{k=1}^n c_{ij}^k e_k,$$
then the algebra $\tilde{\mathfrak g}$ from (\ref{exactseq}) can be defined by  its basis 
$\tilde e_1, \dots, \tilde e_n, \tilde e_{n{+}1}$ 
and structure relations
\begin{equation}
\begin{split}
[\tilde e_i, \tilde e_j]=\sum_{k=1}^n c_{ij}^k \tilde e_k + c(\tilde e_i,\tilde e_j)\tilde e_{n{+}1}, \; 1 \le i < j \le n;\\ 
[\tilde e_i, \tilde e_{n+1}]=0, \; i=1, \dots, n.
\end{split}
\end{equation}

Now let $\tilde{\mathfrak{g}}$ be a filiform Lie algebra, it has one-dimensional center $Z(\tilde{\mathfrak{g}})$
and we have the following one-dimensional central extension:
$$0 \to \mathbb K=Z(\tilde{\mathfrak{g}}) \to \tilde {\mathfrak g}  \to \tilde {\mathfrak g}/ Z(\tilde{\mathfrak{g}})  \to 0 $$
where the quotient Lie algebra  $\tilde {\mathfrak g}/ Z(\tilde{\mathfrak{g}})$ is also filiform Lie algebra.

\begin{proposition}
Let $\mathfrak g$ be a filiform Lie algebra  
 and let its center $Z({\mathfrak{g}})$ be spanned by some $\xi \in Z({\mathfrak{g}})$.
Then $\tilde {\mathfrak g}$ taken from a one-dimensional central extension 
$$0 \to \mathbb K \to \tilde {\mathfrak g}  \to \mathfrak{g} \to 0 $$
with cocycle $c \in \Lambda^2(\mathfrak g)$ 
 is a filiform Lie algebra if and only if
the restricted function $f( \cdot)=c( \cdot, \xi)$ is non-trivial in $\mathfrak{g}^*$.
\end{proposition} 

Let $e_1, \dots, e_n \in \mathfrak g$ be an adapted basis of $\mathfrak{g}$ and
$e^1, \dots, e^n \in {\mathfrak g}^*$ be its dual basis. Then the cocycle $c$ that
corresponds to some filiform central extension can be expressed as
$$ c=(\alpha_1 e^1 + \alpha_2 e^2 ) \wedge e^n + c', \alpha_1^2+ \alpha_2^2 \neq 0, 
c' \in \Lambda^2(e^1,\dots, e^{n{-}1}).$$

So the problems of construction and classification of filiform Lie algebras have an inductive nature:
we start with first dimension $n=3$ and then go ahead computing $H^2(\mathfrak{g})$ and
studing the classes of isomorphisms. 

We restrict now ourselves to the 
special case of $\mathbb N$-graded filiform Lie algebras.     
The corresponding cocycle $c$ will be expressed as
$$ c=\alpha_1 e^1 \wedge e^n + c', \alpha_1 \neq 0, 
c' \in \Lambda^2_{(n{+}1)}(e^1,\dots,e^{n{-}1}).$$

\begin{lemma}
\label{cohomol_m0}
Let $n \ge 3$ then  the cohomology classes of $[ \frac{n+1}{2} ]$ cocycles
$$e^1 \wedge e^n, \:
\left\{
\frac{1}{2}\sum_{i=2}^{2k{-}1} (-1)^i e^i \wedge e^{2k{+}1{-}i}, \:
\; k=2, \dots, \left[ \frac{n{+}1}{2}\right] \right\}.
$$ 
represent a basis of $H^2(\mathfrak m_0(n))$.
\end{lemma}

\begin{corollary}
Up to an isomorphism there are following $\mathbb N$-graded filiform Lie algebras obtained
as one-dimensional central extensions of $\mathfrak m_0(n)$:

a) if $n=2k+1$ -- the algebra $\mathfrak m_0(n+1)$;

b) if $n=2k$ -- two algebras: $\mathfrak m_0(n+1)$
and the algebra $\mathfrak m_{0,1}(2k+1)$ defined by its basis and commutating relations
\begin{equation}
\begin{split}
[e_1,e_i]=e_{i+1}, \quad i=2,\dots,2k{-}1; \\
[e_l,e_{2k-l+1}]=(-1)^{l+1}e_{2k+1}, \quad l=2,\dots,k.
\end{split}
\end{equation}
\end{corollary}

\begin{lemma}[see also \cite{V1}]
\label{cohomol_m2}
$\dim H^2(\mathfrak{m}_2(n))=3$ if $n \ge 5$ and the cohomology classes of cocycles
$$e^1 \wedge e^n+e^2 \wedge e^{n-1}, \quad
e^2 \wedge e^3, \quad
e^2 \wedge e^5-e^3 \wedge e^4$$
represent the basis of $H^2(\mathfrak{m}_2(n))$. 
\end{lemma}

\begin{corollary}
If $n \ge 7$ then up to an isomorphism there exists only one $\mathbb N$-graded filiform Lie algebra obtained
as one-dimensional central extension of $\mathfrak{m}_2(n)$ -- the algebra $\mathfrak{m}_2(n+1)$.
\end{corollary}

\begin{lemma}[see \cite{Mill1}, \cite{Mill2}]
\label{cohomol_Vn}
Let $n \ge 5$ then one can consider a basis in $H^2({\mathcal V}_n)$ consisting of cohomology classes of
following cocycles:
$$
\Omega_{n{+}1}=\frac{1}{2}\sum_{i{+}j{=}n{+}1} (j{-}i) e^i \wedge e^{j},
\quad e^2 \wedge e^3, \quad
e^2 \wedge e^5- 3e^3 \wedge e^4.$$
\end{lemma}

\begin{corollary}
If $n \ge 7$ then up to an isomorphism there is only one $\mathbb N$-graded filiform Lie algebra obtained
as one-dimensional central extension of ${\mathcal V}_n$ and it is isomorphic to the ${\mathcal V}_{n+1}$.
\end{corollary}

Now we start our inductive procedure.

1) the first dimension in which we have a filiform Lie algebra is $n=3$ and there is the only one
algebra -- $3$-dimensional Heisenberg algebra that is isomorphic to $\mathfrak{m}_0(3)$.
 $H^2_{(4)}(\mathfrak{m}_0(3))$
is spanned by $e^1 \wedge e^3$ and we have the only one  central extension --
$\mathfrak{m}_0(4)$.

2) $H^2_{(5)}(\mathfrak{m}_0(4))$ is spanned by two cocycles $e^1 \wedge e^4$ and $e^2 \wedge e^3$ 
and an admissible central extension corresponds to a cocycle of the form $$\gamma e^1 \wedge e^4 + \beta e^2 \wedge e^3, \gamma \ne 0.$$
If $\beta = 0$ we obtain  $\mathfrak{m}_0(5))$ and we'll get the algebra isomorphic to $\mathfrak{m}_2(5))$ if $\beta \ne 0$. 

3) $H^2_{(6)}(\mathfrak{m}_0(5))$ is spanned by $e^1 \wedge e^5$ and  $H^2_{(6)}(\mathfrak{m}_2(5))$ is spanned by 
$e^1 \wedge e^5 + e^2 \wedge e^4$. So there are two corresponding $6$-dimensional filiform Lie algebras:
$\mathfrak{m}_0(6)$ and $\mathfrak{m}_2(6)$.

4) $H^2_{(7)}(\mathfrak{m}_0(6))$ is again spanned by two cocycles $e^1 \wedge e^6$, $e^2 \wedge e^5 - e^3 \wedge e^4$ 
and admissible central extension corresponds to cocycle of the form 
$$\gamma e^1 \wedge e^6 + \beta (e^2 \wedge e^5 - e^3 \wedge e^4), \gamma \ne 0.$$
The case  $\beta = 0$ gives us the algebra $\mathfrak{m}_0(7))$ and if $\beta \ne 0$ we get an algebra that is 
isomorphic to $\mathfrak{m}_{0,1}(7))$.

$H^2_{(7)}(\mathfrak{m}_2(6))$ is spanned by $e^1 \wedge e^6 + e^2 \wedge e^5$ and 
$e^2 \wedge e^5 - e^3 \wedge e^4$. Hence have the following admissible cocycles:
$$\gamma (e^1 \wedge e^6 + e^2 \wedge e^5) + \beta (e^2 \wedge e^5 - e^3 \wedge e^4), \gamma \ne 0.$$
An algebra that corresponds to $\beta = 0$ is isomorphic to $\mathfrak{m}_2(7))$ but if $\beta \ne 0$ taking $\gamma=1$
we obtain an one-parameter family of filiform Lie algebras 
\begin{equation}
\begin{split}
[e_1,e_2]=e_3, \; [e_1,e_3]=e_4, \; [e_1,e_4]=e_5, \\
[e_1,e_5]= e_6, \; [e_1,e_6]=e_7, \; [e_2,e_3]=e_5, \\
[e_2,e_4]=e_6, \; [e_2,e_5]= (1+\beta)e_7, \;
[e_3,e_4]= -\beta e_7.
\end{split}
\end{equation}
Introducing a new parameter $\alpha = - \frac{2\beta+1}{\beta}$ and taking a new basis
$e_1'=e_1$, $e_2'=(2+\alpha)e_2$, $e_3'=(2+\alpha)e_3$, $e_4'=(2+\alpha)e_4$, $e_5'=(2+\alpha)e_5$, $e_6'=(2+\alpha)e_6$,
$e_7'=(2+\alpha)e_7$ we get well-known formulae (see \cite{AG})  for commutating relations of
one-parameter family $\mathfrak{g}_{7, \alpha}$ of non-isomorphic $7$-dimensional Lie algebras:
\begin{equation}
\label{7-family}
\begin{split}
[e_1',e_2']=e_3', \; [e_1',e_3']=e_4', \;
[e_1',e_4']=e_5', \; [e_2',e_3']=(2{+}\alpha) e_5', \;
[e_1',e_5']= e_6',  \\ [e_2',e_4']=(2{+}\alpha) e_6', \;
[e_1',e_6']=e_7', \; [e_2',e_5']=(1{+}\alpha)e_7', \; [e_3',e_4']= e_7'.
\end{split}
\end{equation}

\begin{proposition}
Let us consider an algebra $\mathfrak{g}_{7, \alpha}$ with some $\alpha$.
The subspace $H^2_{(8)}(\mathfrak{g}_{7, \alpha})$ is one-dimensional and it is spanned by
the cohomology class of
$$
e^1 \wedge e^7+
\alpha e^2 \wedge e^6
+e^3 \wedge e^5.
$$
\end{proposition}

Let us denote the corresponding central extension by $\mathfrak{g}_{8, \alpha}$. To get its commutating relations we
have to add to (\ref{7-family}) the following ones:
\begin{equation}
\label{8-fam}
[e_1,e_7]=e_8, \quad [e_2,e_6]=\alpha e_8, \quad [e_3,e_5]= e_8.
\end{equation}

\begin{proposition}
If $\alpha=-\frac{5}{2}$ then $H^2_{(9)}(\mathfrak{g}_{8, \alpha})=0$. If $\alpha \ne -\frac{5}{2}$
then $H^2_{(9)}(\mathfrak{g}_{8, \alpha})$ is one-dimensional and it is spanned by the cohomology class of
$$
e^1 \wedge e^8+
\frac{2\alpha^2+3\alpha-2}{2\alpha+5} e^2 \wedge e^7+
\frac{2\alpha+2}{2\alpha+5}e^3 \wedge e^6+
\frac{3}{2\alpha+5}e^4 \wedge e^5,
$$
\end{proposition}

We denote by $\mathfrak{g}_{9, \alpha}$ the corresponding central extension (for $\alpha \ne -\frac{5}{2}$).

\begin{proposition}
Suppose that $\alpha \ne -\frac{5}{2}$. Then $H^2_{(10)}(\mathfrak{g}_{9, \alpha})$ is one-dimensional and 
it is spanned by the cohomology class of
$$
e^1 \wedge e^9+
\frac{2\alpha^2+\alpha-1}{2\alpha+5} e^2 \wedge e^8+
\frac{2\alpha-1}{2\alpha+5}e^3 \wedge e^7+
\frac{3}{2\alpha+5}e^4 \wedge e^6,
$$
\end{proposition}

We denote by $\mathfrak{g}_{10, \alpha}$ the corresponding central extension.

\begin{proposition}
Let $\alpha   \ne -\frac{5}{2}$. If $\alpha \in \{ -1, -3 \}$  then $H^2_{(11)}(\mathfrak{g}_{10, \alpha})=0$ and if 
$\alpha \ne -1, -3$
then $H^2_{(11)}(\mathfrak{g}_{10, \alpha})$ is one-dimensional and it is spanned by the cohomology class of
$$
\begin{array}{l}
e^1 \wedge e^{10}+
\frac{2\alpha^3{+}2\alpha^2{+}3}{2(\alpha^2{+}4\alpha{+}3)}
e^2 \wedge e^9+
\frac{4\alpha^3{+}8\alpha^2{-}8\alpha{-}21}
{2(\alpha^2{+}4\alpha{+}3)(2\alpha{+}5)}e^3 \wedge e^8+
\\ \qquad \qquad {}+
\frac{3(2\alpha^2{+}4\alpha{+}5)}
{2(\alpha^2{+}4\alpha{+}3)(2\alpha{+}5)}e^4 \wedge e^7+
\frac{3(4\alpha{+}1)}
{2(\alpha^2{+}4\alpha{+}3)(2\alpha{+}5)}e^5 {\wedge} e^6.
\end{array}
$$
\end{proposition}

We denote by $\mathfrak{g}_{11, \alpha}$ the corresponding central extension  for 
 $\alpha \ne -\frac{5}{2}, -1, -3$.

\begin{remark}
The algebras ${\mathcal V}_n, n=7,8,9,10,11$ are isomorphic to $\mathfrak{g}_{7, \alpha_0}$, $\mathfrak{g}_{8, \alpha_0}$, 
$\mathfrak{g}_{9, \alpha_0}$, $\mathfrak{g}_{10, \alpha_0}$, $\mathfrak{g}_{11, \alpha_0}$ respectively with $\alpha_0 = 8$.
We recall also that $\mathfrak{m}_0(3) \cong  \mathfrak{m}_2(3) \cong {\mathcal V}_3$,
$\mathfrak{m}_0(4) \cong  \mathfrak{m}_2(4) \cong {\mathcal V}_4$,
$\mathfrak{m}_2(5) \cong {\mathcal V}_5$, $\mathfrak{m}_2(6) \cong {\mathcal V}_6$.
\end{remark}

\begin{proposition}
Let $\alpha \ne -\frac{5}{2}, -1, -3$ then $H^2_{(12)}(\mathfrak{g}_{11, \alpha}) \ne 0$ if and only if $\alpha=8$. 
$H^2_{(12)}(\mathfrak{g}_{11, 8})$ is one-dimensional and the corresponding
central extension is isomorphic to ${\mathcal V}_{12}$.
\end{proposition}

\begin{proposition}
Besides $H^2_{(n{+}1)}(\mathfrak{g}_{n, \alpha})$ there are only two non-trivial homogeneous subspaces of $H^2(\mathfrak{g}_{n, \alpha})$:
$H^2_{(5)}(\mathfrak{g}_{n, \alpha})$ and $H^2_{(7)}(\mathfrak{g}_{n, \alpha})$ ($n=7-11$). They are spanned respectively by
\begin{equation}
e^2 \wedge e^3, \quad e^2 \wedge e^5 - e^3 \wedge e^4.
\end{equation}
\end{proposition}

Now we have to find out what happends with one-dimensional central extensions of $\mathfrak m_{0,1}(2k+1)$
and of $\mathfrak m_{0,1}(7)$ in particular.

\begin{proposition}
Let $2k+1 \ge 7$ then $H^2_{(2k+2)}(\mathfrak{m}_{0,1}(2k+1))$ is spanned by  the cohomology class of
$$(1-k)e^1 \wedge e^{2k+1} +
\frac{1}{4}\sum_{i{+}j{=}2k{+}2}(-1)^{i{+}1}(j-i) e^{i} \wedge e^{j}.$$
\end{proposition}

We denote by $\mathfrak{m}_{0,2}(2k+2)$ the corresponding central extension. It can be defined
by its basis and commutating relations:
\begin{equation}
\begin{split}
[e_1,e_i]=e_{i+1}, \quad i=2,\dots,2k+1;\\
[e_l,e_{2k-l+1}]=(-1)^{l+1}e_{2k+1}, \quad l=2,\dots,k; \\
[e_j,e_{2k-j+2}]=(-1)^{j+1}(k-j+1)e_{2k+2}, \quad j=2,\dots,k.
\end{split}
\end{equation}

\begin{proposition}
Let $2k+2 \ge 8$ then $H^2_{(2k+3)}(\mathfrak{m}_{0,2}(2k+2))$ is spanned by  
$$(1-k)e^1 \wedge e^{2k+2} +
\frac{1}{4}\sum_{i{+}j{=}2k{+}3}(-1)^{i}(i-2)(j-2) e^{i} \wedge e^{j}.$$
\end{proposition}

We denote by $\mathfrak{m}_{0,3}(2k+3)$ the corresponding central extension. It can be defined
by its basis and commutating relations:
\begin{equation}
\begin{split}
[e_1,e_i]=e_{i{+}1}, \quad i=2,\dots,2k{+}2; \\
[e_l,e_{2k{-}l{-}1}]=({-}1)^{l{+}1}e_{2k{-}1}, \quad l=2,\dots,k{-}1; \\
[e_j,e_{2k{-}j}]=({-}1)^{j{+}1}(k{-}j{-}1)e_{2k}, \quad j=2,\dots,k{-}1; \\
[e_m,e_{2k{-}m{+}1}]=
({-}1)^{m}\left((m{-}2)(k{-}1){-}\frac{(m{-}2)(m{-}1)}{2}\right)e_{2k{+}1},
\quad m=3,\dots,k{-}1.
\end{split}
\end{equation}

\begin{proposition}
Let $2k+3 \ge 9$ then $H^2_{(2k{+}4)}(\mathfrak{m}_{0,3}(2k+3)) = 0$ and therefore 
$\mathfrak{m}_{0,3}(2k+3)$
has no filiform central extensions.
\end{proposition}

We came to the main
\begin{theorem}
\label{osnovn}
Let $\mathfrak{g}=\bigoplus_{i=1}\mathfrak{g}_i$ be
a filiform ${\mathbb N}$-graded Lie, i.e. a graded Lie algebra such that
 \begin{equation}
\begin{split}
\dim \mathfrak{g}_i=1,  \qquad 1 \le i \le n;\\
\mathfrak{g}_i =0, \hspace{11.3mm} \qquad i > n ;\\
[\mathfrak{g}_1,\mathfrak{g}_i]=\mathfrak{g}_{i{+}1}, \quad \hspace{7.8mm}
i \ge 2.
\end{split}
\end{equation}
'Then  $\mathfrak{g}$ is isomorphic to the one and only one Lie algebra from the following list:

1)  Lie algebras of the six sequences $\mathfrak{m}_0(n)$, $\mathfrak{m}_2(n)$,
${\mathcal V}_n$, $\mathfrak{m}_{0,1}(2k{+}2)$,
$\mathfrak{m}_{0,2}(2k{+}2)$, $\mathfrak{m}_{0,3}(2k{+}3)$,
defined by the basis 
$e_1, \dots, e_n$ and commutating relations:
$$
\begin{tabular}{|c|c|c|}
\hline
&&\\[-10pt]
 algebra  & dimension & commutating relations\\
&&\\[-10pt]
\hline
&&\\[-10pt]
$\mathfrak{m}_0(n), \; n \ge 3$ & $n$ &
$[e_1, e_i]=e_{i{+}1}, \quad i=2, \dots, n{-}1$\\
&&\\[-10pt]
\hline
&&\\[-10pt]
$\mathfrak{m}_2(n), \; n \ge 5$ & $n$ &
\begin{tabular}{c}
$[e_1, e_i]=e_{i{+}1}, \quad i=2, \dots, n{-}1$;\\
$[e_2, e_i]=e_{i{+}2}, \quad i=3, \dots, n{-}2$.\\
\end{tabular}
\\
&&\\[-10pt]
\hline
&&\\[-10pt]
${\mathcal V}_n, \; n \ge 12$ & $n$ &
$[e_i,e_j]= \left\{\begin{array}{r}
   (j-i)e_{i{+}j},  i+j \le n;\\
   0, i+j > n.\\[3pt]
   \end{array} \right. $
\\
&&\\[-10pt]
\hline
&&\\[-10pt]
\begin{tabular}{c}
$\mathfrak{m}_{0,1}(2k{+}1)$, \\
$k \ge 3$\\ \end{tabular}
& $2k{+}1$ &
\begin{tabular}{c}
$[e_1, e_i]=e_{i{+}1}, \quad i=2, \dots, 2k$;\\
$[e_l,e_{2k-l+1}]=(-1)^{l+1}e_{2k+1}, \quad l=2,\dots,k$.\\
\end{tabular}
\\
&&\\[-10pt]
\hline
&&\\[-10pt]
\begin{tabular}{c}
$\mathfrak{m}_{0,2}(2k{+}2)$, \\
$k \ge 3$\\ \end{tabular}
& $2k{+}2$ &
\begin{tabular}{c}
$[e_1, e_i]=e_{i{+}1}, \; i{=}2, \dots, 2k{+}1$;\\
$[e_l,e_{2k{-}l{+}1}]=({-}1)^{l{+}1}e_{2k{+}1}, \; l{=}2,\dots,k$;\\
$[e_j,e_{2k{-}j{+}2}]=({-}1)^{j{+}1}(k{-}j{+}1)e_{2k{+}2}, \; j{=}2,
\dots,k$.\\
\end{tabular}
\\
&&\\[-10pt]
\hline
&&\\[-10pt]
\begin{tabular}{c}
$\mathfrak{m}_{0,3}(2k{+}3)$, \\
$k \ge 3$\\ \end{tabular}
& $2k{+}3$ &
\begin{tabular}{c}
$[e_1, e_i]=e_{i{+}1}, \quad i{=}2, \dots, 2k{+}2$;\\
$[e_l,e_{2k{-}l{+}1}]=({-}1)^{l{+}1}e_{2k{+}1}, \; l{=}2,\dots,k$;\\
$[e_j,e_{2k{-}j{+}2}]=({-}1)^{j{+}1}(k{-}j{+}1)e_{2k{+}2}, \; j{=}2,
\dots,k$;\\
$[e_m,e_{2k{-}m{+}3}]=
({-}1)^{m}\left((m{-}2)k-\frac{(m{-}2)(m{-}1)}{2}\right)e_{2k{+}3}$,\\
$m{=}3,\dots,k{+}1$.\\[2pt]
\end{tabular}
\\
\hline
\end{tabular}
$$

2)  Lie algebras of $5$ one-parameter families
$\mathfrak{g}_{n,\alpha}$
of dimensions $n=7-11$ respectively,
defined by their basises and Lie structure relations:
' 
$$
\begin{tabular}{|c|c|c|}
\hline
&&\\[-10pt]
algebra & dimension & commutating relations\\
&&\\[-10pt]
\hline
&&\\[-10pt]
$\mathfrak{g}_{7, \alpha} $ & $7$ &
(1) \\
&&\\[-10pt]
\hline
&&\\[-10pt]
$\mathfrak{g}_{8, \alpha} $ & $8$ &
(1), (2) \\
&&\\[-10pt]
\hline
&&\\[-10pt]
$\mathfrak{g}_{9, \alpha},\: \alpha \ne {-}\frac{5}{2} $ & $9$ &
(1), (2), (3)  \\
&&\\[-10pt]
\hline
&&\\[-10pt]
$\mathfrak{g}_{10, \alpha},\: \alpha \ne {-}\frac{5}{2} $ & $10$ &
(1) -- (4) \\
&&\\[-10pt]
\hline
&&\\[-10pt]
$\mathfrak{g}_{11, \alpha},\, \alpha {\ne} {-}\frac{5}{2},{-}1,{-}3 $ & $11$ &
(1) -- (5) \\[2pt]
\hline
\end{tabular}
$$

' 
$$
\begin{tabular}{|c|c|}
\hline
&\\[-10pt]
& commutating relations for $\mathfrak{g}_{n, \alpha}$\\
&\\[-10pt]
\hline
&\\[-10pt]
(1) & \begin{tabular}{c}
$[e_1,e_2]=e_3, \; [e_1,e_3]=e_4, \;
[e_1,e_4]=e_5, \; [e_1,e_5]= e_6, \; [e_1,e_6]=e_7,$ \\
$[e_2,e_3]=(2{+}\alpha) e_5, \; [e_2,e_4]=(2{+}\alpha) e_6, \;
[e_2,e_5]=(1{+}\alpha)e_7, \; [e_3,e_4]= e_7$\\
\end{tabular}\\
&\\[-10pt]
\hline
&\\[-10pt]
(2) & $[e_1,e_7]=e_8, \quad [e_2,e_6]=\alpha e_8, \quad [e_3,e_5]= e_8$ \\
&\\[-10pt]
\hline
&\\[-10pt]
(3) &
$[e_1,e_8]=e_9, \; [e_2,e_7]=\frac{2\alpha^2{+}3\alpha{-}2}{2\alpha{+}5}e_9, \;
[e_3,e_6]= \frac{2\alpha{+}2}{2\alpha{+}5}e_9, \;
[e_4,e_5]= \frac{3}{2\alpha{+}5}e_9$\\
&\\[-10pt]
\hline
&\\[-10pt]
(4) &
$[e_1,e_9]=e_{10}, \;
[e_2,e_8]=\frac{2\alpha^2{+}\alpha{-}1}{2\alpha{+}5}e_{10}, \;
[e_3,e_7]= \frac{2\alpha{-}1}{2\alpha{+}5}e_{10}, \;
[e_4,e_6]= \frac{3}{2\alpha{+}5}e_{10}$\\
&\\[-10pt]
\hline
&\\[-10pt]
(5) & \begin{tabular}{c}
$[e_1,e_{10}]=e_{11}, \quad
[e_2,e_9]=\frac{2\alpha^3{+}2\alpha^2{+}3}{2(\alpha^2{+}4\alpha{+}3)}e_{11},
\quad [e_3,e_8]= \frac{4\alpha^3{+}8\alpha^2{-}8\alpha{-}21}
{2(\alpha^2{+}4\alpha{+}3)(2\alpha{+}5)}e_{11},$ \\
$[e_4,e_7]=\frac{3(2\alpha^2{+}4\alpha{+}5)}
{2(\alpha^2{+}4\alpha{+}3)(2\alpha{+}5)}e_{11}, \quad
[e_5,e_6]=\frac{3(4\alpha{+}1)}
{2(\alpha^2{+}4\alpha{+}3)(2\alpha{+}5)}e_{11}$\\[2pt]
\end{tabular}\\[2pt]
\hline
\end{tabular}$$
\end{theorem}

\begin{remark}
We can prove Fialowski's  Theorem \ref{fialowski} now as a corollary of the Theorem \ref{osnovn} by considering corresponding 
$\lim {\mathfrak{g}}_n$.
It appears to be useful because \cite{Fial} contains only a sketch of the proof. 
\end{remark}

\begin{remark}
It was proved by Y.~Khakimdjanov in \cite{Kh1} in a different way for $\mathbb K =\mathbb C$,
that there is only a finite number of $\mathbb N$-graded filiform Lie algebras (see also \cite{GozeKh})
in dimensions $\ge 12$, but there were missed three classes $\mathfrak{m}_{0,1}(2k{+}1)$, $\mathfrak{m}_{0,2}(2k{+}2)$,
$\mathfrak{m}_{0,3}(2k{+}3)$. The one-parameter families $\mathfrak{g}_{7, \alpha} $, $\mathfrak{g}_{8, \alpha}$ 
are well-known (see ~\cite{AG} for instance). 
\end{remark}

\section{Symplectic structures and $\mathbb N$-graded filiform Lie algebras}

\begin{lemma}
\label{sympl_hmg}
Let $\mathfrak g$ be a $2k$-dimensional symplectic $\mathbb N$-graded filiform Lie algebra. 
Then its symplectic structure
$\omega$ can be decomposed as $\omega = \omega_{2k{+}1} + \omega'$, where
homogeneous component $\omega_{2k{+}1} \in \Lambda^2_{(2k+1)}(\mathfrak g)$ is 
symplectic form and $\omega'$ is an arbitrary sum of $2$-cocycles of weights
less than $2k+1$.
\end{lemma}

Let us remark that the weight $\lambda$ of the $k$-th power $\omega^k=Ce^1\wedge e^2 \wedge \dots \wedge e^n,\: C\ne 0,$ is equal to 
$\frac{n(n+1)}{2}=k(2k+1)$, but from the other hand   
the maximal possible weight $\lambda_{max}$ of a cocycle in $H^2(\mathfrak g)$ of a $2k$-dimensional 
$\mathbb N$-graded filiform Lie algebra
$\mathfrak g$ is equal to $2k+1$ and hence $\omega^k=\omega_{2k{+}1}^k$ where 
$\omega=\omega_{2k{+}1} + \omega', \; \omega_{2k{+}1} \in H^2_{(2k{+}1)}(\mathfrak g)$.

\begin{theorem}
\label{list_sympl}
We have the following list of non-isomorphic symplectic $\mathbb N$-graded filiform Lie algebras
and their homogeneous symplectic forms (up to non-zero constant):
\begin{equation} 
\label{list_sympl1}
\begin{split}
\mathfrak{m}_0(2k), \; k \ge 2, 
\quad \omega_{2k{+}1}=e^1 {\wedge} e^{2k} + \beta( \sum_{i{=}2}^k ({-}1)^ie^i {\wedge} e^{2k{-}i{+}1} ), \; 
\beta \in \mathbb R, \beta \ne 0;\\
{\mathcal V}_{6}, \: {\mathcal V}_{2k}, \; k \ge 6, \quad  
 \omega_{2k+1}=(2k{-}1) e^1 {\wedge} e^{2k} + (2k{-}3)e^2 {\wedge} e^{2k{-}1} + \dots + e^k {\wedge} e^{k{+}1}; \\
 {\mathfrak g}_{8,\alpha}, \; \alpha \in \mathbb R,
\alpha \ne -\frac{5}{2}, -2, -1, \frac{1}{2}, \\
\omega_{9}(\alpha)=e^1{\wedge}e^8+\frac{2\alpha^2{+}\alpha{-}1}{2\alpha{+}5}e^2{\wedge}e^7+
\frac{2\alpha{-}1}{2\alpha{+}5}e^3{\wedge}e^6+\frac{3}{2\alpha{+}5}e^4{\wedge}e^5;\\
{\mathfrak g}_{10,\alpha}, \; \alpha \in \mathbb R,
\alpha \ne -\frac{5}{2}, -\frac{1}{4}, -1, -3, \alpha_1, \alpha_2, \\
\omega_{11}(\alpha)=e^1 {\wedge} e^{10}+
\frac{2\alpha^3{+}2\alpha^2{+}3}{2(\alpha^2{+}4\alpha{+}3)}e^2 {\wedge} e^9+
\frac{4\alpha^3{+}8\alpha^2{-}8\alpha{-}21}
{2(\alpha^2{+}4\alpha{+}3)(2\alpha{+}5)}e^3 {\wedge} e^8+\\
+\frac{3(2\alpha^2{+}4\alpha{+}5)}
{2(\alpha^2{+}4\alpha{+}3)(2\alpha{+}5)}e^4 {\wedge} e^7+
\frac{3(4\alpha{+}1)}{2(\alpha^2{+}4\alpha{+}3)(2\alpha{+}5)}e^5 {\wedge} e^6.
\end{split}
\end{equation}
where $\alpha_1 \approx -0,893$ and $\alpha_2 \approx 1,534$ are real roots of the polynomials $2\alpha^3+2\alpha+3$ and
$4\alpha^3+8\alpha^2-8\alpha-21$ respectively.
\end{theorem}

The proof follows from the Lemma~\ref{sympl_hmg} and the cohomological data that we used in the proof of the Theorem ~\ref{osnovn}.

\begin{remark}
The classification of symplectic $6$-dimensional nilpotent Lie algebras (based on Morosov's classification)
was done in ~\cite{GozeB} (see also ~\cite{GozeKh}), also in ~\cite{GozeB} it was considered 
the one-parameter family $\mathfrak{g}_{8, \alpha}$ 
and corresponding symplectic form $\omega_{9}(\alpha)$. In \cite{GJKh2}
symplectic (but over $\mathbb C$) low-dimensional ($\dim \mathfrak g \le 10$) filiform Lie algebras were classified .
\end{remark}

Now we can write down formulae (up to an exact form $d \xi$) for an arbitrary symplectic form for all considered above algebras:
\begin{equation}
\label{list_sympl2}
\begin{split}
1) \quad \mathfrak{m}_0(2k): \quad 
\omega= \gamma\omega_{2k{+}1} + \sum_{l{=}2}^{k{-}1} \gamma_{2l{+}1}( \sum_{i{=}2}^l ({-}1)^ie^i {\wedge} e^{2l{-}i{+}1}), 
\; \gamma \ne  0; \\
2) \quad {\mathcal V}_{2k}: \quad  
\omega=\gamma\omega_{2k{+}1}+ \gamma_5e^2 {\wedge} e^{3} + \gamma_7(e^2 {\wedge} e^{5} -3  e^3 {\wedge} e^{4}),
\gamma \ne  0; \\ 
3) \quad  {\mathfrak g}_{8,\alpha}: \quad
\omega(\alpha)=\gamma\omega_{9}(\alpha)+\gamma_5e^2 {\wedge} e^{3}+\gamma_7(e^2 {\wedge} e^{5}-e^3 {\wedge} e^{4}),
\gamma \ne  0; \\
4) \quad {\mathfrak g}_{10,\alpha}: \quad
\omega(\alpha)=\gamma\omega_{11}(\alpha)+\gamma_5e^2 {\wedge} e^{3}+\gamma_7(e^2 {\wedge} e^{5}-e^3 {\wedge} e^{4}),
\gamma \ne  0;
\end{split}
\end{equation}
where $\gamma_5, \gamma_7, \dots, \gamma_{2n{-}3}$ are arbitrary real parameters. 

\section{Symplectic structures and filtered filiform Lie algebras}

\begin{lemma}
Let $\mathfrak g$ be an $2k$-dimensional symplectic filiform Lie algebra, then 
\begin{equation}
\label{nec_cond}
{\rm gr}_C\mathfrak{g} \cong \mathfrak{m}_0(2k).
\end{equation}
\end{lemma}

For the proof we apply the spectral sequence $E$ that corresponds
to the canonical filtration $C$ of $\mathfrak g$.
The term $E_1$ is isomorphic to the cohomology  $H^*({\rm gr}_C \mathfrak g)$.
According to the Theorem \ref{V_ne1} the algebra ${\rm gr}_C \mathfrak g$ is isomorphic either to
$\mathfrak{m}_0(2k)$ or to $\mathfrak{m}_1(2k)$.
But using again the Theorem \ref{V_ne1} we see that 
$\mathfrak m_1(2k)$ has no filiform one-dimensional central extension. 
In other words it has no $2$-cocycles $c$ of maximal filtration
$c=\xi \wedge e^{2k}+c'$ (it can been checked out directly). 
Hence all cocycles in $\Lambda^2(\mathfrak m_1(2k))$
are elements of $\Lambda^*(e^1, \dots, e^{2k{-}1})$. Therefore they correspond  
to degenerate bilinear forms and we have no symplectic cocyle in $H^2(\mathfrak m_1(2k))$. 

It follows from (\ref{nec_cond}) that
$\mathfrak g$ has following commutating relations with respect to some adapted basis $e_1,e_2, \dots, e_{2k}$: 
$$[e_i,e_j]=\sum_{l=0}^{n{-}i{-}j}c^l_{ij}e_{i{+}j{+}l}.$$
And as we remarked above this basis defines a filtration $L$ of ${\mathfrak g}$:
$$L^1=\mathfrak{g} \supset \dots \supset L^p=Span(e_p, \dots, e_{2k}) \supset L^{p{+}1} \supset \dots \supset L^{2k{+}1}=\{0\}.$$

We recall also that $\{c^0_{ij} \}$ are the structure constants of $\mathbb N$-graded filiform Lie algebra
${\rm gr}_L\mathfrak{g}$ and we can consider $\mathfrak{g}$ as a special deformation of
${\rm gr}_L\mathfrak{g}$ that we will call a filtered deformation of a 
$\mathbb N$-graded filiform Lie algebra ${\rm gr}_L\mathfrak{g}$.
But we will see that the condition (\ref{nec_cond}) is not sufficient.       

Let us consider the spectral sequence $E_r$ that corresponds to the filtration $L$. 
The term $E_1$ coincides with the cohomology
of the associated $\mathbb N$-graded Lie algebra ${\rm gr}_L\mathfrak g$:
$$E^{p,q}_1 = H^{p{+}q}_{(-p)}({\rm gr}_L \mathfrak g).$$  

A symplectic cocycle in a symplectic filiform Lie algebra has a maximal filtration and hence we come to the
\begin{theorem}
\label{sympl_def}
Let $\mathfrak g$ be an $2k$-dimensional symplectic filiform Lie algebra. Then $\mathfrak g$
is a filtered deformation of some 
$\mathbb N$-graded symplectic filiform Lie algebra $\mathfrak{g}_0$. 
For corresponding symplectic form $\omega(\mathfrak g)$ we have the following decomposition 
$$\omega(\mathfrak{g})=\omega_{2k{+}1,0} + \sum_{i{+}j{<}2k{+}1} \gamma_{ij}e^i \wedge e^j,$$
where by $\omega_{2k{+}1,0}$ we denote the homogeneous symplectic form of the $\mathfrak{g}_0$
from the list (\ref{list_sympl1}).  
\end{theorem}

\begin{theorem}
\label{criterion}
Let $\mathfrak g$ be a deformation of some $\mathbb N$-graded symplectic filiform Lie algebra
${\mathfrak g}_0$. 
Suppose $e_1,e_2, \dots, e_{2k}$  
is an adapted basis that defines a filtration $L$. By $E_r$ we denote  the corresponding spectral sequence.  
Then the Lie algebra $\mathfrak g$ is symplectic if and only if some homogeneous symplectic
class $[\omega_{2k{+}1}] \in E_1^{{-}2k{-}1, 2k{+}3}=H^{2}_{(2k{+}1)}({\rm gr}_L \mathfrak g)$ survives to the term $E_{\infty}$. 
\end{theorem}

\begin{remark}
If ${\mathfrak g}_0$ is isomorphic
to one of the algebras 
$${\mathcal V}_{2k}, \; {\mathfrak g}_{8, \alpha}, \; {\mathfrak g}_{10, \alpha}$$ 
the space $H^2_{(2k{+}1)}({\mathfrak g}_0)$ is one-dimensional and
the last condition is equivalent to the following one:
all differentials $d_r^{{-}2k{-}1, 2k{+}3}:E_r^{{-}2k{-}1, 2k{+}3} \to E_r^{{-}2k{-}1{+}r, 2k{-}r{+}2}$  are trivial ones.

1) the space $H^{3}({\mathfrak g}_{8, \alpha})$ is generated by homogeneous cocyles of weights $11$, $12$, $13$, $15$;

2) the space $H^{3}({\mathfrak g}_{10, \alpha})$ is generated by homogeneous cocyles of weights $12$, $13$, $14$, $15$;

3) the space $H^{3}({\mathcal V}_{2k})$ is generated by homogeneous cocyles of weights $12$, $15$, $2k{+}3$, $2k{+}4$, $2k{+}5$
(see \cite{Mill2} for details)

Thus we can conclude that the differentials $d_r^{{-}2k{-}1, 2k{+}3}=0$  if 
${\mathfrak g}_0={\mathfrak g}_{8, \alpha}$
or ${\mathfrak g}_{10, \alpha}$. We have only two possibly non-trivial differentials if ${\mathfrak g}_0={\mathcal V}_{2k}$:
$d_{2k{+}1{-}12}^{{-}2k{-}1, 2k{+}3}$ and $d_{2k{+}1{-}15}^{{-}2k{-}1, 2k{+}3}$.

If ${\mathfrak g}_0 \cong \mathfrak{m}_0(2k)$ then the space $H^2_{(2k{+}1)}({\mathfrak g}_0)$ is
two-dimensional, but the cocycle $c=\frac{1}{2}\sum_{i{+}j=2k{+}1} (-1)^i e^i \wedge e^{j}$ 
very rarely survives to the $E_{\infty}$.
\end{remark}

\begin{example}
Let be a filiform Lie algebra $\mathfrak g$ with abelian commutant
$[\mathfrak g, \mathfrak g]$. As it was shown in ~\cite{Br} $\mathfrak g$ has following
commutating relations with respect to some adapted basis $e_1, \dots, e_n$:
\begin{equation}
\begin{split}
[e_1, e_i]=e_{i{+}1}, i=2 ,\dots, n;\\
[e_2, e_j]=e_{j{+}2{+}t}+\sum_{r{=}1} \alpha_r e_{j{+}2{+}t{+}r}, j=3, \dots, n, 
\end{split}
\end{equation}
for some fixed integer $t \ge 0$. 
By changing the basis one can make $\alpha_{t{+}1}=0$.  
An algebra with parameters $\alpha_1, \dots, \alpha_{n{-}5{-}t}$ is isomorphic to another one
with ${\alpha'}_1, \dots, {\alpha'}_{n{-}5{-}t}$ if and only if 
there exists $a \in \mathbb R$ such that ${\alpha'}_i=a^i{\alpha}_i, i=1, \dots, n{-}5{-}t$.

It is obvious that if $t=0$ then $\mathfrak g$ is a filtered deformation of
${\mathfrak m}_2(n)$ and hence $\mathfrak g$ has no symplectic structure for $n=2k \ge 8$.
Now let us consider $t \ge 1$ then  $\mathfrak g$ is a filtered deformation
of another $\mathbb N$-graded filiform Lie algebra -- ${\mathfrak m}_0(n)$. 
\end{example}

\begin{proposition}
Let $\mathfrak g$ be a filiform symplectic Lie algebra  with abelian commutant
$[\mathfrak g, \mathfrak g]$. Then we can choose an adapted basis $e_1, \dots, e_{2k}$
such that: 
\begin{equation}
\begin{split}
[e_1, e_i]=e_{i{+}1}, \:\; i=2 ,\dots, 2k{-}1;\\
[e_2, e_5]=e_{2k}, \: [e_2, e_4]=e_{2k{-}1}+{\alpha}_1 e_{2k},\\
[e_2, e_3]=e_{2k{-}2}+ {\alpha}_1 e_{2k{-}1}+ {\alpha}_2 e_{2k} 
\end{split}
\end{equation}
for some $\alpha_1, \alpha_2 \in \mathbb R$. 
As for symplectic structure, one can take a deformed symplectic structure of $\mathfrak{m}_0(2k)$:
$$\omega=e^1 {\wedge} e^{2k} + \beta( \sum_{i{=}2}^k ({-}1)^ie^i {\wedge} e^{2k{-}i{+}1} )+e^2{\wedge}e^6.$$
If $2k=8$
we can make $\alpha_2=0$ by changing our adapted basis. An
algebra with parameters $\alpha_1, \alpha_2$ is isomorphic to another one
with ${\alpha'}_1, {\alpha'}_2$ if and only if 
there exists $a \in \mathbb R$ such that ${\alpha'}_1=a{\alpha}_1$, ${\alpha'}_2=a^2{\alpha}_2$.
\end{proposition}

\begin{example}
An algebra defined by its basis $e_1, \dots, e_{10}$ and commutating relations
\begin{equation}
\label{pred_primer}
\begin{split}
[e_1, e_i]=e_{i{+}1}, \:\; i=2 ,\dots, 9;\\
[e_2, e_6]=e_{10}, \: [e_2, e_5]=e_{9}+{\alpha}_1 e_{10},\\
[e_2, e_4]=e_{8}+ {\alpha}_1 e_{9}+ {\alpha}_2 e_{10},\\ 
[e_2, e_3]=e_{7}+ {\alpha}_1 e_{8}+ {\alpha}_2 e_{9}+ {\alpha}_3 e_{10}
\end{split}
\end{equation}
does not admit a symplectic structure for all 
$\alpha_1, \alpha_2, \alpha_3 \in \mathbb R$ although it is a deformation
of symplectic $\mathfrak{m}_0(10)$. 
\end{example}
Indeed, $d_1 \equiv 0$ and operator $d_2$ acts on $H^*(\mathfrak{m}_0(10))$ cohomology. In particular
$$d_2^{{-}11,13}([e^2 {\wedge} e^9-e^3 {\wedge} e^8+e^4 {\wedge} e^7-e^5 {\wedge} e^6])={-}2[e^2{\wedge}e^3{\wedge}e^4] \neq 0.$$
Hence as $d_2^{{-}11,13}([e^1{\wedge}e^{10}])=0$ no one symplectic class of $\mathfrak{m}_0(10)$
survives to $E_{\infty}$. 

\section{Contact structures and other applications}

\begin{definition}
A $1$-form $\beta$ on  a $2k{+}1$-dimensional Lie algebra $\mathfrak{g}$
is called contact form if $\beta \wedge (d \beta)^k \ne 0$.
\end{definition}

Let $\mathfrak{g}$ be a nilpotent symplectic Lie algebra, one can consider
a one-dimensional central extension that corresponds to the symplectic cocycle $\omega$ 
$$
0 \to Span(\xi)=\mathbb K \to \tilde {\mathfrak g}  \to \mathfrak{g} \to 0. 
$$

\begin{proposition}
The $1$-form $\beta \in \tilde {\mathfrak g}^*$ such that $\beta( \xi)=1$  
and $\beta ( \mathfrak{g} )=0$ is a contact form on $\tilde {\mathfrak g}$. 
\end{proposition}

\begin{remark}
The previous Proposition is an algebraic version of the Theorems \ref{TO3.2} and \ref{TO3.3}.
\end{remark}

\begin{corollary}
\label{list_cont}
We have the following list of non-isomorphic contact $\mathbb N$-graded filiform Lie algebras:
\begin{equation} 
\label{list_cont1}
\begin{split}
\mathfrak{m}_{0,1}(2k{+}1), \; k \ge 2, \;
{\mathcal V}_{2k{+}1}, \; k \ge 6;  \\ {\mathfrak g}_{7,\alpha}; \quad
 {\mathfrak g}_{9,\alpha}, \; \alpha \in \mathbb R, 
\alpha \ne -\frac{5}{2}, -2, -1, \frac{1}{2}; \\
{\mathfrak g}_{11,\alpha}, \; \alpha \in \mathbb R,
\alpha \ne -\frac{5}{2}, -\frac{1}{4}, -1, -3, \alpha_1, \alpha_2; 
\end{split}
\end{equation}

where $\alpha_1 \approx -0,893$ and $\alpha_2 \approx 1,534$ are real roots of the polynomials $2x^3+2x+3$ and
$4x^3+8x^2-8x-21$ respectively.
\end{corollary}

One can also get versions of main theorems of previous section on filtered algebras but
we want now to discuss possible applications to the problem of classification of
filiform Lie algebras. Namely in odd-dimensional case we have the only one class
of filiform Lie algebras such that ${\rm gr}_C\mathfrak{g}=\mathfrak{m}_0(2k{+}1)$
and two classes in even-dimensional case: ${\rm gr}_C\mathfrak{g}=\mathfrak{m}_0(2k)$
and ${\rm gr}_C\mathfrak{g}=\mathfrak{m}_1(2k)$. Then the next step is to
devide the class ${\rm gr}_C\mathfrak{g}=\mathfrak{m}_0(n)$ into subclasses
that correspond to different ${\rm gr}_L\mathfrak{g}$.   

Another application concerns the examples of nilmanifolds admitting Anosov diffeomorphisms (see \cite{Lr} for references). 

\begin{definition}
A nilpotent Lie algebra $\mathfrak g$ is said to be Anosov if it admits a
hyperbolic automorphism $A$ (i.e. all  eigenvalues of $A$ have absolute value different from $1$) and a basis such that
the matrix of $A$ 
has integer entries and its structure constants are rational numbers (with respect to this basis).  
\end{definition}

Using  Malcev's construction of nilmanifolds we get a nilmanifold $G/\Gamma$ admitting Anosov diffeomorphism
that corresponds to $A$.
 
\begin{theorem}[G.~Lauret \cite{Lr}]
Let $\mathfrak g$ be a $\mathbb N$-graded Lie algebra of finite dimension.
Then the direct sum $\mathfrak{g} \oplus \mathfrak{g}$ is an Anosov Lie algebra.
\end{theorem}

Hence taking $\mathbb N$-graded filiform Lie algebras from the list of the Theorem \ref{osnovn} one can get
new examples of Anosov Lie algebras.

\section*{Acknowledgments}
This work started from discussions on results of \cite{Mill1}, \cite{Mill2}, \cite{Sal} with S.~Salamon
during the visit to the Torino University in December 2000. 
The author is very grateful to A.~Papadopoulos for the invitation to the
UFR de Math\'ematique et d'Informatique de l'Universit\'e Louis Pasteur in Strasbourg where this article was finished.


\begin{thebibliography}{A}

\bibitem{AG}
J.M.~Ancoch\'ea-Bermud\'ez, M.~Goze,
{\it Classification des alg\'ebres des Lie filiformes
de dimension $8$},
Arch. Math., {\bf 50} (1988), 511--525.

\bibitem{BT1}
I.\,K. Babenko and I.\,A.~Taimanov,
{\it On the existence of nonformal simply connected symplectic
manifolds},
Russian Math. Surveys, {\bf 53}:4 (1998), 1082--1083.

\bibitem{BT2}
I.\,K. Babenko and I.\,A.~Taimanov,
{\it On nonformal simply connected symplectic
manifolds}, 
SFB 2888, Preprint 358, 1998.

\bibitem{BG}
C.~Benson, C.~Gordon, 
{\it K\"ahler and symplectic structures on nilmanifolds},
Topology, {\bf 27} (1988), 513--518.

\bibitem{Br}
F.~Bratzlavsky,
{\it Classification des alg\`ebres de Lie nilpotentes de dimension $n$, de classe $n-1$,
dont l'id\'eal d\'eriv\'e est commutatif},
Acad. Roy. Belg. Bull. Cl. Sci. (5) {\bf 60} (1974), 858--865.

\bibitem{Chu}
B.~Chu, 
{\it Symplectic homogeneous spaces},
Trans. Amer. Math. Soc., {\bf 197} (1974), 145--159.

\bibitem{CFG}
L.A.~Cordero, M.~Fernandez, A.~Gray, 
{\it Symplectic manifolds with no K\"aler structure},
Topology, {\bf 25} (1986), 375--380.

\bibitem{Fial}
A.~Fialowski, {\it Classification of graded Lie algebras with two
generators}, Mosc. Univ. Math. Bull., {\bf 38}:2 (1983), 76--79.

\bibitem{GJKh1}
J.R.~G\'omez, A.~Jimen\'ez-Merch\'an, Y.~Khakimdjanov,
{\it Low-dimensional filiform Lie algebras},
J.Pure Appl.Algebra, {\bf 130} (1998), 133--158.

\bibitem{GJKh2}
J.R.~G\'omez, A.~Jimen\'ez-Merch\'an, Y.~Khakimdjanov,
{\it Symplectic structures on the filiform Lie algebras},
J.Pure Appl.Algebra, {\bf 156} (2001), 15--31.

\bibitem{GozeB}
M.~Goze, A.~Bouyakoub,
{\it Sur les alge\`bres de Lie munies d'une forme symplectique},
Rend. Sem. Fac. Sc. Univ. Cagliari, {\bf 57}:1 (1987), 85--97.

\bibitem{GozeKh}
M.~Goze, Y.~Khakimdjanov,
{\it Nilpotent Lie algebras}, Kluwer Academic, Dordrecht, 1996.

\bibitem{Kh1}
 You.B.~Hakimjanov,
{\it Vari\'et\'e des lois d'alg\'ebres de Lie nilpotentes},
Geom. Dedicata, {\bf 40} (1991), 269--295.

\bibitem{Kh2}
 Y.~Khakimdjanov, K.~Khakimdjanova,
{\it Sur une classe d'alg\`ebres de Lie de dimension infinie},
Comm. in Algebra, {\bf 29}:1 (2001), 177--191.

\bibitem{Lr}
G.~Lauret, 
{\it Examples of Anosov diffeomorphisms},
arXiv:math. DS/0202126. 

\bibitem{LM}
A.~Lichnerowicz, A.~Medina, 
{\it On Lie groups with left-invariant symplectic or K\"alerian structures},
Lett. Math. Phys. {\bf 16} (1988), 225--235. 

\bibitem{K}
J.-L. Koszul,
{\it Homologie et cohomologie des alg\`ebres de Lie},
Bull. Soc. Math. France {\bf 78} (1950), 65--127.

\bibitem{Mal}
A. Malcev,
\textit{On a class of homogeneous spaces},
Amer. Math. Soc. Transl. (1) {\bf 9} (1962), 276--307.

\bibitem{Mill1}
D.V. Millionschikov,
\textit{Cohomology of nilmanifolds and Gontcharova's theorem},
Russian Math. Surveys, {\bf 56}:4 (2001), 153--154 
(in russian). 

\bibitem{Mill2}
D.V. Millionschikov,
\textit{Cohomology of nilmanifolds and Gontcharova's theorem},
in "Global Differential geometry: The Mathematical Legacy
of Alfred Gray", M. Fernandez and J.Wolf ed., AMS CONM 288 (2001),
381--385.

\bibitem{Mill3}
D.V. Millionschikov,
\textit{$\mathbb N$-graded filiform Lie algebras},
Russian Math. Surveys, {\bf 57}:3 (2002), (to appear in russian). 

\bibitem{Mor}
V.~Morosov,
{\it Classification of nilpotent Lie algebras of order $6$},
Izv. Vyssh. Uchebn. Zaved. Mat., {\bf 4} (1958), 161--171.

\bibitem{Nz}
K. Nomizu,
\textit{On the cohomology of homogeneous spaces of nilpotent Lie groups},
Ann. of Math. {\bf 59} (1954), 531--538.

\bibitem{Sal}
S.M.~Salamon, 
{\it Complex structures on nilpotent Lie algebras},
J. Pure Appl. Algebra  {\bf 157} (2001), 311--333. 

\bibitem{TO}
A.~Tralle, J.~Oprea, 
{\it Symplectic manifolds with no K\"aler structure},
Lect. Notes in Math., 1661, Springer, 1997. 

\bibitem{V1}
M.~Vergne,
{\it R\'eductibilit\'e de la vari\'et\'e des alg\`ebres de Lie nilpotentes},
C.R. Acad. Sc. Paris. {\bf 263} (1966), 4--6.

\bibitem{V2}
M.~Vergne,
{\it Cohomologie des alg\`ebres de Lie nilpotentes},
Bull. Soc. Math. France {\bf 98} (1970), 81--116.

\end{thebibliography}
\end{document}